\documentstyle[12pt]{article}
\begin{document}
\bibliographystyle{alpha}
\baselineskip=5.4mm

\begin{center}
{\large\bf  Biclique Coverings and the Chromatic Number}
\ \\
\ \\
\begin{tabular}{ccc}
{\large Dhruv Mubayi} &   {\large Sundar Vishwanathan} \\

Department of Mathematics, &      Department of Computer Science \\
Statistics, and Computer Science &      Indian Institute of Technology\\
University of Illinois &      Mumbai\\  Chicago, IL 60607, USA &
India 400076 \\

{\tt mubayi@math.uic.edu}      & {\tt sundar@cse.iitb.ernet.in} \\
\end{tabular}
\end{center}

\def\Box{{\rm \mbox{\,\rule[1.5pt]{.1pt}{4pt}\rule[1.5pt]{4pt}{.1pt}\hskip -4pt\rule[5.5pt]{4pt}{.1pt}\rule[1.5pt]{.1pt}{4pt}}}\,\,}
\newtheorem{defi}{Definition}
\newenvironment{definition}{\begin{defi}\rm}{\end{defi}}
\newtheorem{claim}{Claim}

\newtheorem{theorem}{Theorem}
\newtheorem{lemma}[theorem]{Lemma}
\newtheorem{fact}[theorem]{Fact}
\newtheorem{proposition}[theorem]{Proposition}
\newtheorem{note}[theorem]{Note}
\newtheorem{corollary}[theorem]{Corollary}
\newtheorem{example}[theorem]{Example}
\newtheorem{conjecture}[theorem]{Conjecture}
\newcommand{\blackslug}{\hbox{\hskip 1pt \vrule width 4pt height 8pt
depth 1.5pt \hskip 1pt}}
\newcommand{\QED}{\quad\blackslug\lower 8.5pt\null}
\newenvironment{proof}{\noindent {\bf Proof:}}{\QED}
\newenvironment{sketch}{\noindent {\bf Proof (sketch):}}{\QED}
\newcommand{\st}{\,|\;}
\newcommand{\inter}{\cap }
\newcommand{\union}{\cup }
\newcommand{\edge}{\relbar\joinrel\joinrel\relbar}
\newcommand{\qed}{\quad\blackslug\lower 8.5pt\null}
\newcommand{\nospaceqed}{\blackslug\lower 8.5pt\null}
\newcommand{\floor}[1]{\left\lfloor #1 \right\rfloor}
\newcommand{\ceil}[1]{\left\lceil  #1 \right\rceil}


\begin{abstract}
Consider a graph  $G$  with chromatic number $k$ and  a
collection of complete bipartite graphs, or bicliques,
that cover the edges of $G$. We prove
the following two results:
\medskip

\noindent
 $\bullet$ If the bicliques   partition the
edges of $G$, then their number is at least $2^{\sqrt{\log_2 k}}$.
This is the first improvement of the easy lower bound of $\log_2 k$,
while the Alon-Saks-Seymour conjecture states that this can be
improved to $k-1$.
\medskip

\noindent $\bullet$ The sum of the orders of the bicliques is at
least $(1-o(1))k\log_2 k$. This generalizes, in asymptotic form, a
result of Katona and Szemer\'edi who proved that the minimum is
$k\log_2 k$ when $G$ is a clique.

\end{abstract}

\section{Introduction}
It is a well-known fact that the minimum number of bipartite graphs
needed to cover the edges of a graph $G$ is $\lceil \log \chi(G)
\rceil$, where $\chi(G)$ is the chromatic number of $G$ (all logs
are to the base 2).  Call a complete bipartite graph  a {\em
biclique}. Two classical theorems study related questions. One is
the Graham-Pollak theorem \cite{GP} which states that the minimum
number of bicliques needed to partition $E(K_k)$ is $k-1$. Another
is the Katona-Szemer\'edi theorem \cite{KS}, which states that the
minimum of the sum of the orders of a collection of bicliques that
cover $E(K_k)$ is $k\log k$. Both of these results are best
possible.

An obvious way to generalize these theorems is to ask whether the
same results hold for any $G$ with chromatic number $k$.

\begin{conjecture} \label{assc}{\bf (Alon-Saks-Seymour)}
The minimum number of bicliques needed to partition the edge set of
a graph $G$ with chromatic number $k$ is $k-1$.
\end{conjecture}
Note that every graph has  a partition of this size, simply by
taking a proper coloring $V_1, \ldots V_k$ and letting the $i$th
bipartite graph be $(V_i, \cup_{j>i} V_j)$.

Another motivation for  Conjecture \ref{assc} is that the
non-bipartite analogue is an old conjecture of Erd\H
os-Faber-Lov\'asz.  The Erd\H os-Faber-Lov\'asz conjecture remains
open although it has been proved asymptotically by Kahn~\cite{K}.
Conjecture \ref{assc} seems much harder than the Erd\H
os-Faber-Lov\'asz conjecture, indeed, as far as we know there are no
nontrivial results towards it except the folklore lower bound of
$\log_2 k$ which doesn't even use the fact that we have a partition.
Our first result improves this to a superlogarithmic bound for $k$
large.

\begin{theorem} \label{ass}
The number of bicliques needed to partition the edge set of a graph
$G$ with chromatic number $k$ is at least $2^{\sqrt{2\log
k}(1+o(1))}$.
\end{theorem}

Motivated by Conjecture \ref{assc}, we make the following conjecture
that generalizes the Katona-Szemer\'edi theorem.

\begin{conjecture} \label{ksconj}
Let $G$ be a graph with chromatic number $k$. The  sum of the orders
of any collection of bicliques that cover the edge set of $G$  is at
least $k\log k$.
\end{conjecture}

We prove Conjecture \ref{ksconj} asymptotically.

\begin{theorem} \label{ks}
 Let $G$ be a graph with chromatic number $k$, where $k$ is sufficiently
large. The  sum of the orders of any collection of bicliques that
cover the edge set of $G$  is at least
$$k\log k-k\log\log k -k \log \log \log k.$$
\end{theorem}

The next two sections contain the proofs of Theorems \ref{ass} and
\ref{ks}.

\section{The Alon-Saks-Seymour Conjecture}
It is more convenient to phrase and prove our result in inverse
form.
 Let $G$ be a disjoint union of
$m$ bicliques $(A_i, B_i), 1\leq i \leq m$. The Alon-Saks-Seymour
conjecture then states that the chromatic number of $G$ is at most
$m+1$.

We prove the following theorem which immediately implies Theorem
\ref{ass}.
\begin{theorem}
Let $G$ be a disjoint union of $m$ bicliques. Then $\chi(G) \leq
m^{\frac{1+\log m}{2}} (1+o(1))$.
\end{theorem}
{\bf Proof.} We will begin with a proof of a worse bound. We will
first show that $\chi(G) \leq  m^{\log m} (1+o(1))$. A color will be
an ordered tuple of length at most $\log m$,  with each element a
positive integer of value at most  $m$. We will construct this tuple
in stages. In the $i$th stage we will fill in the $i$th co-ordinate.
Note that the length of the tuple may vary with vertices.

With each vertex $v$, at stage $i$, we will associate a set
$S(i,v)\subset V(G)$. The set  $S(i,v)$ will contain all vertices
which have the same color sequence, so far, as $v$ (in particular,
$v \in S(i,v)$ for all $i$).

A biclique $(A_j, B_j)$ is said to {\em cut} a subset of vertices
$S$ if $S \cap A_j \not= \emptyset$ and $S \cap B_j \not=
\emptyset$.

Consider two bicliques $(A_k, B_k)$ and $(A_l, B_l)$ from our
collection. Since they are edge disjoint, $(A_l,B_l)$ cuts either
$A_k$ or $B_k$, but not both.

Fix a vertex $v$. We set $S(0,v):=V(G)$. The assignment for the
$i+1$st stage is as follows. Suppose we have defined $S(i,v)$. Let
${\cal F}(i,v)$ denote the set of all bicliques that cut $S(i,v)$.
For each biclique $(A_j, B_j) \in {\cal F}(i,v)$ for which $v \in
A_j \cup B_j$, let $C_j$ be the set among $A_j, B_j$ that contains
$v$ and let $D_j$ be the set among $A_j, B_j$ that omits $v$. For a
vertex $v$, check if there is a biclique $(A_j,B_j) \in {\cal
F}(i,v)$ such that $v \in A_j \cup B_j$ and
\begin{itemize}
  \item The number of bicliques in ${\cal F}(i,v)$ that cut $C_j$ is
     at most the number that cut $D_j$.
\end{itemize}

If  there is such a $j$, then the $i+1$st  co-ordinate of the color
of $v$ is $j$ and $S(i+1,v) = S(i,v) \cap C_j$. If there are many
candidates for $j$, pick one arbitrarily.

If there is no such $(A_j, B_j)$, then  the coloring of $v$ ceases
and the vertex will not be considered in subsequent stages. In other
words, the final color of vertex $v$ will be a sequence of length
$i$.

Note that in this process every vertex is assigned
a color except vertices that were not assigned a color in
the very first step.
We will show below that if a vertex is assigned a color
then this coloring is proper. The same argument shows that
the vertices that do not get assigned a color in the first
step form an independent set.
These vertices  are all assigned a
special color which is swallowed up in the $o(1)$ term.

The following technical lemma establishes the statements needed to
prove correctness and a bound on the number of colors used.

\begin{lemma}
For each vertex $v$, the set $S(i,v)$ is determined by the color
sequence $x_1,\ldots,x_i$ assigned to  the vertex $v$. Also, the
number of bicliques that cut $S(i,v)$ is at most $m/2^i$.
\end{lemma}
{\bf Proof.} The proof is by induction on $i$. Both statements are
trivially true  for $i=0$. For the inductive step,  assume that
$S(i-1,v)$ is determined by $x_1,\ldots,x_{i-1}$ and at most
$m/2^{i-1}$ bicliques cut $S(i-1,v)$. If $v$ ceases to be colored
then we are done by induction. Now suppose that $v$ is colored with
$x_i=t$ in step $i$. Then $(A_t, B_t) \in {\cal F}(i,v)$ and $v \in
A_t \cup B_t$. As before, define $C_t$ and $D_t$. Because $v$ is
colored in this step, the number of bicliques in ${\cal F}(i,v)$
that cut $C_{t}$ is at most the number which cut $D_{t}$. As $S(i,v)
= C_{t} \cap S(i-1,v)$, we see that $S(i,v)$ is determined by $x_1,
\ldots, x_{i-1}, t$. Also, since the number of bicliques that cut
$C_{t}$ is at most half the number that cut $S(i-1,v)$ the second
assertion follows. \qed

We argue first that the coloring is proper. Assume for a
contradiction that two adjacent vertices $v$ and $w$ are assigned
the same color sequence. Suppose the sequence is  of length $i$.
Then by the previous lemma $S(i,v) = S(i,w)$. There has to be one
biclique, say $(A_p,B_p)$, such that $v \in A_p$ and $w\in B_p$. If
the number of bicliques in ${\cal F}(i,v)$ that cut $A_p$ is at most
the number that cut $B_p$ then $v$ will be given a color in the
$i+1$st step. Otherwise $w$ will be colored. In any case, at least
one of them will be given a color contradicting our assumption that
both sequences are of length $i$. This argument also shows that the
vertices which were not assigned a color in the first step form an
independent set.  The coloring stops when ${\cal F}(i,v)$ is empty
for every vertex and that happens after $\log m$ steps from the
lemma.

A simple observation helps in reducing this bound by a square-root
factor. At each stage, the colorings of the $S(i,v)$s are
independent. Hence the colors only matter within the vertices in
each of these sets. The number of bicliques that cut $S(i,v)$ is at
most $m/2^i$. We {\em renumber} these bicliques from $1$ to $m/2^i$.
Hence the labels in the $i$th stage will be restricted to this set.
The total number of colors used, of length $i$ is at most
$m\cdot\frac{m}{2}\cdots\frac{m}{2^{i}}$. The number for $i<m$ is
swallowed up in the $o(1)$ term and the value for $i=m$ simplifies
to the main term in the bound given. \qed

\section{Generalizing the Katona-Szemer\'edi Theorem}

In this section we prove Theorem \ref{ks}. Given a graph $G$, let
$b(G)$ denote the minimum, over all collections of bicliques that
cover the edges of $G$, of the sum of the orders of these bipartite
graphs.

One proof of the Katona-Szemer\'edi theorem is due to
Hansel~\cite{H} and the same proof yields the following lemma which
is part of folklore.

\begin{lemma} Let $G=(V,E)$ be an $n$ vertex graph with independence
number $\alpha$.  Then $b(G) \ge n \log (n/\alpha)$.
\end{lemma}

In other words for any graph $G$, $\alpha(G) \geq
\frac{n}{2^{b(G)/n}}$. Let $k=\chi(G)$. We may assume that $n\leq
k\log k$, since we are done otherwise. Let $G=G_0$. Starting with
$G_0$, repeatedly remove independent sets of size given by Hansel's
lemma as long as the number of vertices is at least $k$. Let the
graphs we get be $G_0, G_1,\ldots, G_t$. Let $|V(G_i)| = n_i$ and
$\beta = \max_i 2^{{b(G_i)/n_i}}$. Let this maximum be achieved for
$i=p$. From the definition, we see that $n_{i+1} \leq n_i
(1-\frac{1}{2^{b(G_i)/n_i}})$. Hence $n_t \leq n
(1-1/\beta)^t<ne^{-t/\beta}<n2^{-t/\beta}$ and  together with
$n_t\ge k$ we obtain
$$t\leq  \beta\log (n/k).$$

There are two cases to consider. First suppose that  $ t \geq k/\log
k$. Then from the above two inequalities we obtain
$$2^{b(G_p)/n_p} \log (n/k) \ge k/\log k.$$
Taking $\log$s and using the facts that $n\leq k\log k$ and $n_p \ge
k$  we get
$$b(G_p) \ge k(\log k- \log\log k-\log\log\log k).$$

We now consider the case that $t < k/\log k$. Let $G'$ be the graph
obtained after removing an independent set from $G_t$. By definition
of $t$ we have $|V(G')| < k$. Also $\chi(G') \geq k(1-1/\log k)$.
Since the color classes of size one in an optimal coloring form a
clique, this implies that $G'$ has a clique of size at least
$k(1-2/\log k)$. Using the fact that $\log(1-x)>-2x$ for $x$
sufficiently small and applying the Katona-Szemer\'edi theorem, we
get $b(G') \ge k\log k -3k$. \qed

Note that in the proof $b(G_i)$ could use different
covers, but with sizes smaller than the one induced by $b(G_0)$.
One can get better lower order terms by adjusting the
threshold between the two cases.

\section{Acknowledgments}

The research of Dhruv Mubayi was supported in part by  NSF grant DMS
0653946.

\end{document}